\documentclass[12pt]{article}

\usepackage{amsmath, graphicx}
\usepackage{amstext,amsfonts,amsbsy,eucal,amssymb, amsthm,color}
\usepackage{tikz, here}
\usetikzlibrary{matrix,graphs,arrows,positioning,calc,decorations.markings,shapes.symbols}

\advance \topmargin by -\headheight
\advance \topmargin by -\headsep
     
\evensidemargin \oddsidemargin
\newtheorem{thm}{Thorem}[section]
\newtheorem{prp}{Proposition}[section]

\theoremstyle{definition}
\newtheorem{definition}{Definition}[section]  

\newtheorem{rem}[definition]{Remark}

\def\Z{\operatorname{\mathbb{Z}}}

\def\C{\operatorname{\mathbb{C}}}
\def\P{\operatorname{\mathbb{P}}}

\def\cX{\operatorname{\mathcal{X}}}
\def\cY{\operatorname{\mathcal{Y}}}

\def\cH{\operatorname{\mathcal{H}}}
\def\cE{\operatorname{\mathcal{E}}}
\def\cK{\operatorname{\mathcal{K}}}

\begin{document}
\begin{center}
{\large {\bf Space of initial conditions for the four-dimensional Garnier system revisited}}\\

\medskip
{\large Tomoyuki Takenawa}

\medskip
{\small
Faculty of Marine Technology, Tokyo University of Marine Science and Technology\\ 2-1-6 Etchu-jima, Koto-ku, Tokyo, 135-8533, Japan\\
E-mail: takenawa@kaiyodai.ac.jp}
\end{center}

\medskip

\begin{abstract}
A geometric study  is given for the 4-dimensional Garnier system.
By the resolution of indeterminacy, the group of its B\"aklund transformations is lifted to a group of pseudo-isomorphisms between rational varieties obtained from $\P^2 \times \P^2$ by 10 or 21 blow-ups. The root basis is discussed in the N\'eron-Severi bilattices for the space with 10 blow-ups.
\end{abstract}

\section{Introduction}
\subsection{Background and overview}
The method of Okamoto-Sakai's spaces of initial conditions (SIC in short) is a  powerful tool for studying  the classification and symmetry of the Painlev\'e equations in two dimensions \cite{Okamoto1979, Sakai2001}.
However, studies of higher dimensional initial value spaces are scarce due to the theoretical and computational complexity. 

The Garnier system is a natural generalization of the Painlev\'e equations to higher order and  obtained as the monodromy preserving deformation for a Fuchsian ODE of the Schlesinger type: 
\begin{align}\label{Fuchsian}
\frac{d}{dx}{\bf y}(x)=\sum_{i=1}^{n+2} \frac{A_i}{x-z_i}{\bf y}(x),\end{align}
where $n$ is a positive integer, ${\bf y}$ is two-dimensional and the eigenvalues of $A_i$'s and $A_\infty=-\sum_{i=1}^{n+2} A_i$ are different with each other. More concretely, the monodromy preserving deformation is given by so called the Schlesinger equation:
\begin{align}
\frac{\partial A_i}{\partial z_j}=\frac{[A_j,A_i]}{z_j-z_i}\quad  (j\neq i), \quad  
\frac{\partial A_i}{\partial z_i}=-\sum_{j\neq i}\frac{[A_j,A_i]}{z_j-z_i} 
\end{align}
[Jimbo-Miwa-M\^ori-Sato 1980].
Since the singular points of the Fuchsian equation $(z_1,$ $\dots,$ $z_{n+2},$ $\infty)$ can be normalized to $(0,1,\infty, z_1',\dots, z_n')$ by a M\"obius transformation for $x$, the Garnier system has $n$ independent variables, or more precisely, they are $n$-compatible systems of ODEs. 
On the other hand the Garnier system has $2n$ dependent variables, since the dimension of the moduli space of ODEs of the form of \eqref{Fuchsian} is $2n$.
It is well known that  the Garnier system has Painlev\'e property, i.e. movable singularities are at most poles \cite{IKSY1991}.

In particular, when $n=2$, the Garnier system is a commuting pair of four-dimensional systems of ordinary differential equations. In \cite{KO1984} H. Kimura and K. Okamoto  showed that it can be written in a polynomial Hamiltonian form of two directions as
$$\frac{dq_i}{ds_j}=\frac{\partial H_j}{\partial p_i},\ 
\frac{dp_i}{ds_j}=-\frac{\partial H_j}{\partial q_i}\quad (i,j=1,2)$$
with the Hamiltonians
\begin{align*}
s_1 (s_1 - 1)H_1=&\Big(q_1 (q_1 - 1) (q_1 - s_1) - 
      \frac{s_1 (s_1 - 1)}{s_1 - s_2} q_1 q_2\Big) p_1^2 \\&+ 
   2 q_1 q_2 \Big(q_1 + \frac{s_1 (s_2 - 1)}{s_1 - s_2}\Big) p_1 p_2 + 
   q_1 q_2 \Big(q_2 - \frac{s_2 (s_1 - 1)}{s_1 - s_2}\Big) p_2^2\\& -
   \Big\{(\kappa_0 - d) q_1 (q_1 - 1) + \kappa_1 q_1 (q_1 - s_1) + \theta_1 (q_1 - 1) (q_1 - s_1) \\&\quad + 
      \theta_2 q_1 \Big(q_1 + \frac{s_1 (s_2 - 1)}{s_1 - s_2}\Big) - 
      \theta_1 \frac{s_1 (s_1 - 1)}{s_1 - s_2} q_2\Big\} p_1\\&
 + \Big((2\alpha_0+\kappa_\infty) q_1 q_2 + 
      \theta_2 q_1 \frac{s_2 (s_1 - 1)}{s_1 - s_2} - \theta_1 q_2 \frac{s_1 (s_2 - 1)}{s_1 - s_2}\Big) p_2 \\& +  \alpha_0(\alpha_0+\kappa_\infty) q_1\\
H_2=& \{\mbox{replacing as $q_1\leftrightarrow q_2$,  $p_1\leftrightarrow p_2$,  $s_1\leftrightarrow s_2$, $\theta_1\leftrightarrow \theta_2$ in $H_1$}\},
\end{align*}
where $\kappa_0$, $\kappa_1$, $\kappa_\infty$, $\theta_1$, $\theta_2$ and $\alpha_0$ are parameters independent from the Hamiltonian flows, and $d$ is given by  $d=2\alpha_0+\kappa_0+\kappa_1+\kappa_\infty+\theta_1+\theta_2$.

The Hamiltonians can be written using the Hamiltonian for the sixth Painlev\'e equation
\begin{align*}
&s(s-1)H_{{\rm VI}}(q,p,s,\kappa_0,\kappa_1,\kappa_\infty,\theta,\alpha_0)\\
=& q(q-1)(q-s)p^2-\Big\{-(2\alpha_0+\kappa_1+\kappa_\infty+\theta)q(q-1)+\kappa_1q(q-s)\\&\qquad +\theta_1(q-1)(q-s)\Big\}p+ \alpha_0(\alpha_0+\kappa_\infty) q
\end{align*}
as
\begin{align*}
s_1(s_1-1)H_1
=& s_1(s_1-1)H_{{\rm VI}}(q_1,p_1,s_1,\kappa_0,\kappa_1,\kappa_\infty,\theta_1,\alpha_0)\\&
+(2q_1p_1+q_2p_2+2\alpha_0+\kappa_\infty )q_1q_2p_2
+\frac{s_1(s_1-1)}{s_1-s_2}(-q_1p_1+\theta_1)q_2p_1\\&
+\frac{s_1(s_2-1)}{s_1-s_2}(2q_1p_1-\theta_1)q_2p_2
+\frac{s_2(s_1-1)}{s_1-s_2}\big(-q_2p_2^2 +\theta_2 (p_2-p_1)\big)q_1.
\end{align*}

In this paper we construct the space of initial conditions for this four-dimensional Garnier system,  where a fiber space $\pi: X\to B$ is called the space of initial conditions for the ODE system $\varphi$ if 
\begin{enumerate}
\item[(i)] $\varphi$ is regularly defined at any point in $X$;
\item[(ii)] For any point $x$ in $X$, any path $\gamma$ in $B$ passing through $\pi(x)$ can be lifted to $ \gamma' \subset X$ by $\varphi$, where $\pi:\gamma' \to \gamma$ is an isomorphism.   
\end{enumerate} 
The most typical examples of the SIC are those for the Painlev\'e equations\cite{Okamoto1979}, while H. Kimura constructed the SIC for the Garnier system in $n$ variables in \cite{Kimura1993}.

In particular, for the sixth Painlev\'e equation (the Garnier system with $n=1$), the SIC was constructed by blowing up a Hirzebrugh surface 8 times and excluding 5 irreducible curves (called vertical leaves),  where $B$ is $\{t \in \C ~|~ t \neq 0,1\}$. 
 For the Garnier system with $n=2$, it was constructed by blowing up a $4$-dimensional minimal projective variety 10 times, where $B$ is $\{(t_1, t_2) \in \C^2 ~|~ t_i \neq 0,1\ (i=1,2),\ t_1\neq t_2 \}$.
In the $n=2$ case, M. Suzuki also constructed by gluing 13 affine spaces $\C^4 $ using the B\"acklund  transformations in \cite{Suzuki2006} and Y. Sasano (reported as he) constructed the SIC by 13 times blowing up from $\P^4$ in \cite{Sasano2007}.  
Note that these systems have the Painlev\'e property, which guarantees that any path in $B$ can be globally lifted.

Our approach to constructing the SIC is similar to Kimura's, but follows the way of H. Sakai constructing the SIC for discrete Painlev\'e systems \cite{Sakai2001}\footnote{A concrete expression of the SIC for the sixth Painlev\'e equation from $\P^1\times \P^1$ can be found in \cite{DST2013} and \cite{KNY2017}, for example.}.   
Our approach differs from Kimura's at the following points:
\begin{enumerate}
\item[(i)] we blow up from $\P^2\times \P^2$  instead of twisted four-dimensional variety;
\item[(ii)] we resolve the indeterminacy of B\"acklund transformations instead of the system of ODEs.
\end{enumerate} 
In the discrete setting, the SIC is a compact variety without excluding vertical leaves.
Thus, not only does this approach provide a good computation prospect, but also
allows us to easily see the homology and cohomology structure and to find the root lattice in this structure.  

For a discrete dynamical system $\varphi_n: X_n\to X_{n+1}$, $n\in \Z$, the sequence of manifolds $\{X_n\}_{n\in \Z}$ is called a SIC for $\{\varphi_n\}_{n\in \Z}$ if $\varphi_n: X_n\to X_{n+1}$ is a pseudo-isomorphism for all $n\in \Z$ (see Section 1.3 below).\footnote{Some of the contents of this paper have already been published in the proceedings for a special lecture of the infinite integrable systems session of the Mathematical Society of Japan in 2021, but this paper is the first publication of the proofs, details and accompanying results.}

\subsection{B\"acklund transformations}

We use birational symmetries (B\"acklund transformations) of the Garnier sytem found by H. Kimura \cite{Kimura1990} and T. Tsuda \cite{Tsuda2003}. Let us denote $(\kappa_0, \kappa_1, \kappa_\infty, \theta_1, \theta_2, \alpha_0, s_1, s_2)$ by $\alpha$ and $w(f)$ by $\bar{f}$ for a birational action $w$. B\"acklund transformations change parameters as Table~\ref{tab:parameters}.

\begin{table}[H]
\begin{center}
 \caption{Action on parameters}
 \label{tab:parameters}

\begin{tabular}{|c|c|c|c|c|c|c|c|c|}\hline
& $\bar{\kappa}_0$&$\bar{\kappa}_1$&$\bar{\kappa}_\infty$&$\bar{\theta}_1$&$\bar{\theta}_2$&$\bar{\alpha}_0$&$\bar{s}_1$&$\bar{s}_2$\\\hline 
$w_{\kappa_0}$&$-\kappa_0$&$\kappa_1$&$\kappa_\infty$&$\theta_1$&$\theta_2$&$\alpha_0+\kappa_0$&$s_1$&$s_2$\\$
w_{\kappa_1}$&$\kappa_0$&$-\kappa_1$&$\kappa_\infty$&$\theta_1$&$\theta_2$&$\alpha_0+\kappa_1$&$s_1$&$s_2$\\$
w_{\kappa_\infty}$&$\kappa_0$&$\kappa_1$&$-\kappa_\infty$&$\theta_1$&$\theta_2$&$\alpha_0+\kappa_\infty$&$s_1$&$s_2$\\$
w_{\theta_1}$&$\kappa_0$&$\kappa_1$&$\kappa_\infty$&$-\theta_1$&$\theta_2$&$\alpha_0+\theta_1$&$s_1$&$s_2$\\$
w_{\theta_2}$&$\kappa_0$&$\kappa_1$&$\kappa_\infty$&$\theta_1$&$-\theta_2$&$\alpha_0+\theta_2$&$s_1$&$s_2$\\$
w_{\alpha_0}$&$d-\kappa_0$&$d-\kappa_1$&$-\kappa_\infty$&$-\theta_1$&$-\theta_2$&$-\alpha_0$&$s_1$&$s_2$\\
$
\sigma_1$&$\kappa_1$&$\kappa_0$&$\kappa_\infty$&$\theta_1$&$\theta_2$&$\alpha_0$&$s_1^{-1}$&$s_2^{-1}$\\
$
\sigma_2$&$\kappa_0$&$\kappa_\infty$&$\kappa_1$&$\theta_1$&$\theta_2$&$\alpha_0$&$\frac{s_1}{s_1-1}$&$\frac{s_2}{s_2-1}$\\$
\sigma_3$&$\kappa_0$&$\kappa_1$&$\theta_1$&$\kappa_\infty$&$\theta_2$&$\alpha_0$&$s_1^{-1}$&$s_1^{-1}s_2$\\
$
\sigma_4$&$\kappa_0$&$\kappa_1$&$\kappa_\infty$&$\theta_2$&$\theta_1$&$\alpha_0$&$s_2$&$s_1$\\ \hline
\end{tabular}
\end{center}
\end{table}
The actions on dependent variables are as Table~\ref{tab:qp},
where
\begin{align}
d&=2\alpha_0+\kappa_0+\kappa_1+\kappa_\infty +\theta_1+\theta_2\nonumber\\
Q_{12}&=q_1+q_2-1\label{Q12}\\
Q_{12}^s&=q_1/s_1+q_2/s_2-1 \label{Q12s}\\ 
P_{12}&=q_1p_1+q_2p_2+\alpha_0. \nonumber
\end{align}

\begin{table}[H]
\begin{center}
 \caption{Action on $(q,p)$ variables}
 \label{tab:qp}

\begin{tabular}{|c|c|c|c|c|}\hline
&$\bar{q}_1$&$\bar{q}_2$&$\bar{p}_1$&$\bar{p}_2$\\\hline
$w_{\kappa_0}$&$q_1$&$q_2$&$p_1-\frac{\kappa_0}{s_1Q_{12}^s}$&$p_2-\frac{\kappa_0}{s_2Q_{12}^s}$\\
$w_{\kappa_1}$&$q_1$&$q_2$&$p_1-\frac{\kappa_1}{Q_{12}}$&$p_2-\frac{\kappa_1}{Q_{12}}$\\
$w_{\kappa_\infty}$&$q_1$&$q_2$&$p_1$&$p_2$\\
$w_{\theta_1}$&$q_1$&$q_2$&$p_1-\theta_1/q_1$&$p_2$\\
$w_{\theta_2}$&$q_1$&$q_2$&$p_1$&$p_2-\theta_2/q_2$\\
$w_{\alpha_0}$&$\frac{s_1p_1(q_1p_1-\theta_1)}{P_{12}(P_{12}+\kappa_\infty)}$ &$\frac{s_2p_2(q_2p_2-\theta_2)}{P_{12}(P_{12}+\kappa_\infty)}$&$-q_1p_1/\bar{q}_1$ &$-q_2p_2/\bar{q}_2$\\ 
$\sigma_1$&$s_1^{-1}q_1$&$s_2^{-1}q_2$&$s_1p_1$&$s_2p_2$\\
$\sigma_2$&$\frac{q_1}{Q_{12}}$&$\frac{q_2}{Q_{12}}$&
$Q_{12}(p_1-P_{12})$&$Q_{12}(p_2-P_{12})$\\
$\sigma_3$&$q_1^{-1}$&$-q_1^{-1}q_2$&$-q_1P_{12}$&$-q_1p_2$\\
$\sigma_4$&$q_2$&$q_1$&$p_2$&$p_1$\\\hline
\end{tabular}
\end{center}
\end{table}

Similar to the case of the sixth Painlev\'e equation, let us introduce new coordinates $(q_i,r_i)=(q_i, q_ip_i)$, then these actions are written more simply as Table~\ref{tab:qr}, where 
\begin{align}
R_{12}&=r_1+r_2+\alpha_0.\label{R12}
\end{align}

\begin{table}[H]
\begin{center}
 \caption{Action on $(q,r)$ variables}
 \label{tab:qr}

\begin{tabular}{|c|c|c|c|c|}\hline
&$\bar{q}_1$&$\bar{q}_2$&$\bar{r}_1$&$\bar{r}_2$\\\hline
$w_{\kappa_0}$&$q_1$&$q_2$&$r_1-\frac{\kappa_0 q_1}{s_1Q_{12}^s}$&$r_2-\frac{\kappa_0 q_2}{s_2Q_{12}^s}$\\
$w_{\kappa_1}$&$q_1$&$q_2$&$r_1-\frac{\kappa_1 q_1}{Q_{12}}$&$r_2-\frac{\kappa_1q_2}{Q_{12}}$\\
$w_{\kappa_\infty}$&$q_1$&$q_2$&$r_1$&$r_2$\\
$w_{\theta_1}$&$q_1$&$q_2$&$r_1-\theta_1$&$r_2$\\
$w_{\theta_2}$&$q_1$&$q_2$&$r_1$&$r_2-\theta_2$\\
$w_{\alpha_0}$&$\frac{s_1r_1(r_1-\theta_1)}{q_1R_{12}(R_{12}+\kappa_\infty)}$ &$\frac{s_2r_2(r_2-\theta_2)}{q_2R_{12}(R_{12}+\kappa_\infty)}$&$-r_1$ &$-r_2$ \\ $\sigma_1$&$s_1^{-1}q_1$&$s_2^{-1}q_2$&$r_1$&$r_2$\\
$\sigma_2$&$\frac{q_1}{Q_{12}}$&$\frac{q_2}{Q_{12}}$&
$r_1- q_1 R_{12}$&$r_2-q_2R_{12}$\\
$\sigma_3$&$q_1^{-1}$&$-q_1^{-1}q_2$&$-R_{12}$&$r_2$\\
$\sigma_4$&$q_2$&$q_1$&$r_2$&$r_1$\\\hline
\end{tabular}
\end{center}
\end{table}

In the next section we will consider $\P^2 \times \P^2$ using this coordinate system.

\subsection{Basic facts}
In this paper, we use the following basic facts about birational maps between higher dimensional varieties; see \S~2 of \cite{CT2019} for details. \\
 
\noindent {\it Pseudo-isomorphisms}\\
Let $\cX$ and $\cY$ be smooth projective varieties.  For a birational map $f:\cX \to \cY$, let $I(f)$ denote the indeterminate set (i.e. the set of points where $f$ is not defined) of $f$ in $\cX$.
 
We say a sequence of birational maps $\varphi_n: \cX_n \to\cX_{n+1}$ for smooth projective varieties $\cX_n$ ($n\in \Z$)   
to be algebraically stable if $$\left(\varphi_{n+k-1}  \circ \cdots \circ  \varphi_{n+1}  \circ \varphi_n\right)^* =  \varphi_n^* \circ \varphi_{n+1}^* \circ \cdots \circ \varphi_{n+k-1}^* $$ 
holds as a mapping from the Picard group of $\cX_{n+k}$ to that of $\cX_n$
for any integers $n$ and $k\geq 1$.

\begin{prp}[\cite{BK2008, Bayraktar2012}] \label{AS2}
A sequence of birational maps $\varphi_n: \cX_n \to\cX_{n+1}$ for smooth projective varieties $\cX_n$ ($n\in {\mathbb Z}$)   
is algebraically stable if and only if
there do not exist integers $n$ and $k\geq 1$ and a divisor $D$ on $\cX_{n-1}$ such that $\varphi(D\setminus I(\varphi_{n-1} ))\subset I(\varphi_{n+k-1}  \circ \cdots \circ  \varphi_{n+1}  \circ \varphi_n)$.\footnote{
This statement is a non-autonomous analog of a proposition shown in \cite{BK2008, Bayraktar2012}. The proof does not change except in notations.}
\end{prp}

We call a birational mapping $f:\cX\to \cY$ a {\it pseudo-isomorphism}
if $f$ is isomorphic except on finite number of subvarieties of codimension two at least. This condition is equivalent to that there is no prime divisor pulled back to the zero divisor by $f$ or $f^{-1}$. Hence, if $\varphi_n$ is a pseudo-isomorphism for each $n$, then $\{\varphi_n\}_{n\in \Z}$ and $\{\varphi_n^{-1}\}_{n\in \Z}$ are algebraically stable.  

\begin{prp}[\cite{DO1988}]\label{NSauto}
Let $\cX$ and $\cY$ be smooth projective varieties and
$\varphi$ a pseudo-isomorphism from $\cX$ to $\cY$. Then
$\varphi$ acts on the N\'eron-Severi bi-lattice as an automorphism preserving the intersections.
\end{prp}
The N\'eron-Severi bi-lattice of a smooth rational variety $\cX$ is isomorphic to $H^2(\cX, \Z) \times H_2(\cX, \Z)$ which is explicitly given in the following.\\

\noindent {\it Blow-ups}\\
Recall that in local coordinates $U\subset \C^N$, the blow-up along a subvariety $V$ of dimension $N-d$, $d\geq 2$, written as 
$$x_1-f_1(x_{d+1},\dots, x_N)=\dots=x_d-f_d(x_{d+1},\dots, x_N)=0,$$
where $f_i$'s are holomorphic functions,  
is a birational morphism $\pi:X\to U$ such that $X=\{U_i\}$ is an open variety given by
\begin{align*}
U_i=\{(u_1^{(i)},\dots, u_d^{(i)},x_{d+1},\dots, x_N )\in \C^N\}\quad (i=1,\dots,d)
\end{align*}
with $\pi: U_i \to U$: 
\begin{align*}(x_1,\dots,x_N)=&(u_1^{(i)}u_i^{(i)}+f_1, \dots, u_{i-1}^{(i)}u_i^{(i)}+f_{i-1},
u_i^{(i)}+f_i, \\
\quad &u_{i+1}^{(i)}u_i^{(i)}+f_{i+1}, \dots,  u_d^{(i)}u_i^{(i)}+f_d, x_{d+1},\dots,x_N).
\end{align*} 
It is convenient to write the coordinates of $U_i$ as
$$\left(\frac{x_1-f_1}{x_i-f_i},\dots,\frac{x_{i-1}-f_{i-1}}{x_i-f_i},x_i-f_i,
\frac{x_{i+1}-f_{i+1}}{x_i-f_i},\dots,\frac{x_d-f_d}{x_i-f_i},x_{d+1},\dots, x_N \right).$$
The exceptional divisor $E$ is written as $u_i=0$ in $U_i$ and each point in the center of  the blow-up corresponds to a subvariety isomorphic to $\P^{d-1}$:
$(x_1-f_1: \dots: x_d-f_d)$. Thus $E$ is locally a direct product $V \times \P^{d-1}$.
We called such $\P^{d-1}$ a fiber of the exceptional divisor. 
(In algebraic setting the affine charts often need to be embedded into higher dimensional space.)\\

\noindent {\it Case of $\P^2 \times \P^2$ }\\
Let $\cX$ be a rational variety obtained by $K$ successive blowups from $\P^2 \times \P^2$  and 
$$({\bf x}_1,{\bf x}_2)=(x_{10}:x_{11}:x_{12},~ x_{20}:x_{21}:x_{22})$$
the direct product of homogeneous coordinate chart. 
Let $\cH_i$ denote the total transform of the class of a hyper-plane ${\bf c}_i{\bf x}_i=c_{i0}x_{i0}+c_{i1}x_{i1}+c_{i2}x_{i2}=0$, where ${\bf c}_i=(c_{i0}:c_{i1}:c_{i2})$ is a constant vector in $\P^{2}$, and $\cE_k$ the total transform of the $k$-th exceptional divisor class.

Let $h_i$ denote the total transforms of the class of a line
$$\{{\bf x}~|~ {\bf x}_j={\bf c}_j (\forall j\neq i),\ {\bf x}_i=s{\bf a}_i+t{\bf b}_i (\exists (s:t)\in \P^1) \},$$ 
where ${\bf a}_i$,  ${\bf b}_i$
and ${\bf c}_j$'s are constant vectors in $\P^2$, and $e_k$ the class of a line in a fiber of the $k$-th blow-up.
Then, the Picard group $\simeq H^2(\mathcal{X},\Z)$ and its Poincar\'e dual $\simeq H_2(\mathcal{X},\Z)$ are lattices
\begin{align}\label{NSbasis}
H^2(\mathcal{X},\Z)=\bigoplus_{i=1}^2 \Z \cH_i \oplus \bigoplus_{k=1}^{K} \Z \cE_k,\quad 
H_2(\mathcal{X},\Z)=\bigoplus_{i=1}^2 \Z h_i \oplus \bigoplus_{k=1}^{K} \Z e_k
\end{align}
and the intersection form is given by
\begin{align}
\langle \cH_i, h_j\rangle = \delta_{ij},\quad 
\langle \cE_k, e_l\rangle = -\delta_{kl},\quad
\langle \cH_i, e_k\rangle =
\langle \cE_k, h_i\rangle =0
\end{align} 
for $i,j=1,2$ and $1\leq k,l\leq K$.

Moreover, the anti-canonical divisor class of $\cX$ is given by
\begin{align}
-\cK_{\cX}= 3 \cH_1 + 3\cH_2 - \sum_{k=1}^K  (3-d_k) \cE_k,
\end{align}
where $d_k$ is the dimension of the center manifold for the $k$-th blow-up.

\section{Construction of the space of initial conditions}

First of all, let us compactify the phase space $(q_1,q_2,p_1,p_2)\in \C^4$ to
$(Q_0:Q_1:Q_2)\times (R_0:R_1:R_2)\in  \P^2\times \P^2$, where the original coordinates correspond to $Q_0\neq 0$ and $R_0\neq 0$ through $(q_1,q_2,r_1,r_2)=(Q_1/Q_0,Q_2/Q_0,R_1/R_0,R_2/R_0)$ as usual.
Let $\cX_\alpha$ be a rational projective variety obtained by successive 10 blow-ups from $\P^2\times \P^2$ with the center of each blow-up $C_i$ ($i=1,\dots,10$)
given as follows. 
\begin{align*}
\begin{array}{ll}
C_1: q_1=r_1=0  &U_1: (u_1,q_2,v_1,r_2)=(q_1,q_2,r_1/q_1,r_2)\\
C_2: q_1=r_1-\theta_1 =0&U_2: (u_2,q_2,v_2,r_2)=(q_1,q_2,(r_1-\theta_1)/q_1,r_2)\\
C_3: q_2=r_2=0  &U_3: (q_1,u_3,r_1,v_3)=(q_1,q_2,r_1,r_2/q_2)\\
C_4:  q_2=r_2-\theta_2 =0&U_4: (q_1,u_4,r_1,v_4)=(q_1,q_2,r_1,(r_2-\theta_2)/q_2)\\
C_{5}:  Q_0=R_{12}=0  &U_5: (u_5,q_{2,1},r_1,v_5)=(1/q_1,q_2/q_1,r_1,q_1R_{12})\\
C_{6}:Q_0=R_{12}+\kappa_\infty=0  &U_6: (u_6,q_{2,1},r_1,v_6)=(1/q_1,q_2/q_1,r_1,q_1(R_{12}+\kappa_\infty))\\
C_{7}:  R_0=Q_{12}=A_{12}=0  &U_{7}: (q_1,u_7, v_7,w_7)=
(q_1,1/r_1, Q_{12}r_1,A_{12}r_1)\\
C_{8}: u_7=v_7-\kappa_1 q_1=0  &U_{8}:  (q_1,u_8, v_8,w_7)=
(q_1,1/r_1, (v_7-\kappa_1 q_1)r_1,w_7)\\
C_{9}:  R_0=Q_{12}^s=A_{12}^s=0  &U_{9}: (q_1,u_9, v_9,w_9)=
(q_1,1/r_1, Q_{12}^s r_1, A_{12}^s r_1)\\
C_{10}: u_9=v_9-\kappa_0 q_1/s_1=0  &U_{10}:  (q_1,u_{10}, v_{10},w_9)=
(q_1, 1/r_1, (v_9-\kappa_0 q_1/s_1)r_1,w_9)
\end{array}
\end{align*}
where $Q_{12}$, $Q_{12}^s$ are $R_{12}$ are \eqref{Q12}, \eqref{Q12s} and \eqref{R12} respectively, and 
\begin{align}
A_{12}&=q_2/q_1-r_2/r_1\\ 
A_{12}^s&= (s_1q_2)/(s_2q_1)-r_2/r_1.
\end{align}
Here, $C_i$ is two-dimensional subvariety for $i\neq7,9$, while $C_7$ and $C_9$ are 1-dimensional.

The following proposition follows immediately from the basic facts in Section 1.3.

\begin{prp}
The Picard group ($\simeq H^2(\cX_\alpha,\Z)$) and its dual $\simeq H_2(\cX_\alpha,\Z)$ are 12-dimensional lattices
\begin{align*}
H^2(\cX_\alpha,\Z)= \Z \cH_q \oplus \Z \cH_r \oplus \bigoplus_{k=1}^{10} \Z \cE_k,\quad 
H_2(\cX_\alpha,\Z)= \Z h_q \oplus \Z h_r \oplus \bigoplus_{k=1}^{10} \Z e_k,
\end{align*}
where $\cH_q$, $\cH_r$ and $\cE_k$ denote the total transforms of the classes for a hyper-plane $c_0Q_0+c_1Q_1+c_2Q_2=0$, $c_0R_0+c_1R_1+c_2R_2=0$ with $(c_0:c_1:c_2)\in \P^2$ and the $k$-th exceptional divisor, and $h_q$, $h_p$ and $e_k$ denote the total transforms of the classes for a generic line 
$$\{ (Q, R)=({\bf c}, s {\bf a} + t {\bf b}) \in \P^2 \times \P^2 ~|~(s:t) \in \P^1\},$$
$$\{ (Q, R)=(s {\bf a} + t {\bf b}, {\bf c}) \in \P^2 \times \P^2 ~|~(s:t) \in \P^1\},$$
and a generic fiber of the $k$-th blow-up. The intersection form is given by
\begin{align}
\langle \cH_i, h_j\rangle = \delta_{ij},\quad 
\langle \cE_k, e_l\rangle = -\delta_{kl},\quad
\langle \cH_i, e_k\rangle=\langle \cE_k, h_i\rangle =0.
\end{align}
The anti-canonical divisor class is 
$$-\cK_{\cX} = 3 \cH_q + 3 \cH_r - \sum_{k=1,2,3,4,5,6,8,10} \cE_k  -2 \cE_7 -2 \cE_9.$$
\end{prp}

Let $A$ denote the set of generic values of the parameter $\alpha$.
The following theorem holds.

\begin{thm}
The B\"acklund transformations $w=w_i$ for $i=\kappa_0,\kappa_,\kappa_\infty,\theta_1,\theta_2$ and  $w=\sigma_j$ for $j=1,2,3,4$ can be lifted to pseudo-isomorphisms from $\mathcal{X}_\alpha$ to $\mathcal{X}_{w(\alpha)}$, where $w$ acts a bijection on $A$. 
\end{thm}

Table~\ref{tab:pic10} is the actions of the above $w$, where we omit the preserved elements and 
$\cE_{i_1,i_2,\dots,i_k}$ is the abbreviation for $\cE_{i_1}+\cE_{i_2}+\cdots+\cE_{i_k}$.

\begin{table}[ht]
\begin{center}
 \caption{Action on the Picard group}
 \label{tab:pic10}

\begin{tabular}{|c|c|}\hline
$w_{\kappa_0}$& $\cH_r\leftrightarrow \cH_q+\cH_r-\cE_{9,10}$\\&
$\cE_9 \leftrightarrow \cH_q-\cE_{10}$, \quad
$\cE_{10} \leftrightarrow \cH_q-\cE_{9}$\\\hline
$w_{\kappa_1}$& $\cH_r\leftrightarrow \cH_q+\cH_r-\cE_{7,8}$\\&
$\cE_7 \leftrightarrow \cH_q-\cE_{8}$, \quad
$\cE_8 \leftrightarrow \cH_q-\cE_{7}$\\\hline
$w_{\kappa_\infty}$& $\cE_5\leftrightarrow \cE_6$\\\hline
$w_{\theta_1}$& $\cE_1\leftrightarrow \cE_2$\\\hline
$w_{\theta_2}$& $\cE_3\leftrightarrow \cE_4$\\\hline
$\sigma_1$& $\cE_7\leftrightarrow \cE_9$, \quad $\cE_8\leftrightarrow \cE_{10}$\\\hline
$\sigma_2$& $\cH_r\leftrightarrow \cH_q+\cH_r-\cE_{5,7}$,\quad
$\cE_5\leftrightarrow \cH_q-\cE_{7}$\\
&$ \cE_6\leftrightarrow \cE_8$,\quad
$\cE_7\leftrightarrow \cH_q-\cE_{5}$ \\\hline
$\sigma_3$& $\cE_1\leftrightarrow \cE_5, \quad \cE_2\leftrightarrow \cE_{6}$\\\hline
$\sigma_4$& $\cE_1\leftrightarrow \cE_3, \quad \cE_2\leftrightarrow \cE_{4}$\\\hline
\end{tabular}
\end{center}
\end{table}

Furthermore, direct calculations also allow us to verify the following theorem.

\begin{thm}
Let $\cX_\alpha^o$ be the open variety obtained by excluding the following proper transforms from $\mathcal{X}_\alpha$:
\begin{align*}
\begin{array}{ll}
Q_1=0:& \cH_q-\cE_1-\cE_2\\
Q_2=0:& \cH_q-\cE_3-\cE_4\\
Q_0=0:& \cH_q-\cE_5-\cE_6\\
R_0=0:& \cH_r-\cE_7-\cE_9\\
R_0=Q_{12}=A_{12}=0:& \cE_7-\cE_8\\
R_0=Q_{12}^s=A_{12}^s=0:&\cE_9- \cE_{10},
\end{array}
\end{align*}
then the family $\{\mathcal{X}_\alpha^o\}_{\alpha \in A}$ is a SIC for the Garnier system.
\end{thm}

However, we need further blow-ups for $w_{\alpha_0}$ which is lifted to a pseudo-isomorphism. Actually, $w_{\alpha_0}$ maps a hypersurface $Q_0=0$ to a subvariety $Q_1=Q_2=0$ whose codimension is two, which means that we need blow-up along $Q_1=Q_2=0$.   
Note that both $Q_0=0$ and $Q_1=Q_2=0$ are included in $\cX_\alpha \setminus \cX_\alpha^o$, vertical leaves. 

Let us blow-up $\P^2\times \P^2$ along $C_{11},\dots,C_{21}$, where
\begin{align*}
\begin{array}{ll}
C_{11}:& Q_0=Q_2=0  \quad U_{11}: (u_{11},v_{11},r_1,r_2)=(1/q_1,q_2,r_1,r_2)\\
C_{12}:& Q_0=Q_1=0  \quad U_{12}: (u_{12},v_{12},r_1,r_2)=(1/q_2,q_1,r_1,r_2)\\
C_{13}:& q_1=q_2=0  \quad U_{13}: (u_{13},v_{13},r_1,r_2)=(q_1,q_2/q_1,r_1,r_2)\\
C_{14}:&  q_1=q_2-1=R_0=R_1=0\\&U_{14}: (x_{14},w_{14},v_{14},u_{14})=(q_1r_2,(q_2-1)r_2,r_1,1/r_2)\\
C_{15}:&  q_1=q_2/s_2-1=R_0=R_1=0\\ &U_{15}:(x_{15},w_{15},v_{15},u_{15})=(q_1r_2,(q_2/s_2-1)r_2,r_1,1/r_2)\\
C_{16}:&  q_2=q_1-1=R_0=R_2=0\\ &U_{16}:(x_{16},w_{16},v_{16},u_{16})=(q_2r_1,(q_1-1)r_1,r_2,1/r_1)\\
C_{17}:&  q_2=q_1/s_1-1=R_0=R_2=0\\&U_{17}:(x_{17},w_{17},v_{17},u_{17})=(q_2r_1,(q_1/s_1-1)r_1,r_2,1/r_1)\\
C_{18}:&  Q_0=Q_1+Q_2=R_0=R_1+R_2=0\\&U_{18}: (x_{18},w_{18},v_{18},u_{18})=(r_1/q_1,(q_2/q_1+1)r_1,r_2+r_1,1/r_1)\\
C_{19}:&  Q_0=Q_1/s_1+Q_2/s_2=R_0=R_1+R_2=0\\&U_{19}: (x_{19},w_{19},v_{19},u_{19})=(r_1/q_1,(q_2/q_1+s_2/s_1)r_1,r_2+r_1,1/r_1)\\
C_{20}: &q_1+s_1 (s_2-1)/(s_1-s_2)=q_2 + s_2 (s_1 - 1)/(s_2-s_1)=R_0 =0\\
&U_{20}: (v_{20},w_{20},u_{20},r_{2,1})\\
 &=((q_1+s_1 (s_2-1)/(s_1-s_2))r_1,(q_2 + s_2 (s_1 - 1)/(s_2-s_1))r_1,1/r_1,r_2/r_1)
\end{array}
\end{align*}
and $C_{21}=w_{\alpha_0}(C_{20})$.

Then, all the B\"acklund transformations in Table~\ref{tab:qr} are lifted to pseudo-isomorphisms as Table~\ref{tab:pic}.

\begin{table}[H]
\begin{center}
 \caption{Action on the Picard group with 21 blow-ups}
 \label{tab:pic}

\begin{tabular}{|c|c|}\hline
$w_{\kappa_0}$& $\cH_r\leftrightarrow \cH_q+\cH_r-\cE_{9,10,15,17,19,20}$\\&
$\cE_9 \leftrightarrow \cH_q-\cE_{10,15,17,19,20}$, \quad
$\cE_{10} \leftrightarrow \cH_q-\cE_{9,15,17,19,20}$\\\hline
$w_{\kappa_1}$& $\cH_r\leftrightarrow \cH_q+\cH_r-\cE_{7,8,14,16,18,20}$\\&
$\cE_7 \leftrightarrow \cH_q-\cE_{8,14,16,18,20}$, \quad
$\cE_8 \leftrightarrow \cH_q-\cE_{7,14,16,18,20}$\\\hline
$w_{\kappa_\infty}$& $\cE_5\leftrightarrow \cE_6$\\\hline
$w_{\theta_1}$& $\cE_1\leftrightarrow \cE_2$\\\hline
$w_{\theta_2}$& $\cE_3\leftrightarrow \cE_4$\\\hline
$w_{\alpha_0}$& $\cH_q\leftrightarrow 2\cH_q+2\cH_r-\cE_{1,2,3,4,5,6,11,12,13,14,15,16,17,18,19}$\\&
$\cE_1\leftrightarrow \cH_r-\cE_{1,14,15},\quad 
\cE_2\leftrightarrow \cH_r-\cE_{2,14,15},\quad 
\cE_3\leftrightarrow \cH_r-\cE_{3,16,17}$\\&
$\cE_4\leftrightarrow \cH_r-\cE_{4,16,17},\quad 
\cE_5\leftrightarrow \cH_r-\cE_{5,18,19},\quad 
\cE_6\leftrightarrow \cH_r-\cE_{6,18,19}$\\&
$\cE_7\leftrightarrow \cE_9, \quad \cE_8\leftrightarrow \cE_{10}$\\
&$\cE_{11}\leftrightarrow \cH_q-\cE_{1,2,12,13,14,15},\quad 
\cE_{12}\leftrightarrow \cH_q-\cE_{3,4,11,13,16,17}$\\& 
$\cE_{13}\leftrightarrow \cH_q-\cE_{5,6,11,12,18,19}\quad
\cE_{20}\leftrightarrow \cE_{21}$\\ \hline
$\sigma_1$& $\cE_7\leftrightarrow \cE_9, \quad \cE_8\leftrightarrow \cE_{10}$\\& $\cE_{14}\leftrightarrow \cE_{15}, \quad \cE_{16}\leftrightarrow \cE_{17},
\quad \cE_{18}\leftrightarrow \cE_{19}$\\\hline
$\sigma_2$& $\cH_r\leftrightarrow \cH_q+\cH_r-\cE_{5,7,14,16,18,19},\quad
\cE_5\leftrightarrow \cH_q-\cE_{7,14,16,18,20}$\\
&$ \cE_6\leftrightarrow \cE_8,\quad
\cE_7\leftrightarrow \cH_q-\cE_{5,11,12,18,19} \quad
\cE_{11}\leftrightarrow \cE_{16}, \quad \cE_{12}\leftrightarrow \cE_{14},
\quad \cE_{19}\leftrightarrow \cE_{20}$\\\hline
$\sigma_3$& $\cE_1\leftrightarrow \cE_5, \quad \cE_2\leftrightarrow \cE_{6}$\\& 
$\cE_{11}\leftrightarrow \cE_{13},\quad \cE_{14}\leftrightarrow \cE_{18}, \quad \cE_{15}\leftrightarrow \cE_{19}$\\\hline
$\sigma_4$& $\cE_1\leftrightarrow \cE_3, \quad \cE_2\leftrightarrow \cE_{4}$\\& 
$\cE_{11}\leftrightarrow \cE_{12},\quad \cE_{14}\leftrightarrow \cE_{16}, \quad \cE_{15}\leftrightarrow \cE_{17}$\\\hline

\end{tabular}
\end{center}
\end{table}

\section{Root system}
Since the space with 21 blow-ups is too complicated to see the structure,  in this section we consider the action on the N\'eron-Severi bi-lattice for the case of 10 blow-ups. 

Let $\cX_{\alpha}$ be the space of initial conditions obtained by the first 10 blow-ups.
Define the root vectors $\alpha_i$ ($i=0,1,\dots,5$) and the co-root vectors as
\begin{align*}
\begin{array}{ll}
\alpha_0=\frac{1}{2}(\cH_q+2\cH_r-2\cE_1-2\cE_3-2\cE_5),&
\alpha_1=\cH_q-\cE_9-\cE_{10}\\
\alpha_2=\cH_q-\cE_7-\cE_8,&
\alpha_3=\cE_5-\cE_6\\
\alpha_4=\cE_1-\cE_2,&
\alpha_5=\cE_3-\cE_4\\
\\
\check{\alpha}_0=h_q-e_1-e_3-e_5,&
\check{\alpha}_1=h_r-e_9-e_{10}\\
\check{\alpha}_2=h_r-e_7-e_8,&
\check{\alpha}_3=e_5-e_6\\
\check{\alpha}_4=e_1-e_2,&
\check{\alpha}_5=e_3-e_4
\end{array}
\end{align*}
(see Figure \ref{fig:Dynkin}).

\begin{figure}[ht]
\begin{center}
   \begin{tikzpicture}[
			elt/.style={circle,draw=black!100,thick, inner sep=0pt,minimum size=2mm}, scale=0.7]
		\path 	(0,0) 	node 	(g0) [elt] {} 
				(0,2) 	node 	(g1) [elt] {}
		        (-1.9,0.62) 	node 	(g2) [elt] {}
		        ( -1.18,-1.92) 	node  	(g3) [elt] {}
		        (1.18,-1.92) 	node  	(g4) [elt] {}
		        ( 1.9,0.62) 		node  	(g5) [elt] {};
		\draw [black,line width=1pt ] (g0) -- (g1) (g0) -- (g2)  (g0) -- (g3) (g0) -- (g4) (g0) -- (g5); 
			\node at ($(g0.east) + (0.3,-0.1)$) 	{$\alpha_{0}$};
			\node at ($(g0.west) + (-0.55,-0.2)$) 	{{\color{red} $-\frac{5}{2}$}};
			\node at ($(g1.north) + (0,0.2)$) 	{$\alpha_{1}$};
			\node at ($(g1.west) + (-0.35,0)$) 	{{\color{red} $-2$}};
			\node at ($(g1.west) + (-0.2,-1)$) 	{{\color{red} $1$}};
			\node at ($(g2.west) + (-0.3,0)$) 	{$\alpha_{2}$};
			\node at ($(g2.north) + (0,0.3)$) 	{{\color{red} $-2$}};
			\node at ($(g2.west) + (0.95,-0.01)$) 	{{\color{red} $1$}};
			\node at ($(g3.west) + (-0.3,0)$) 	{$\alpha_{3}$};
			\node at ($(g3.east) + (0.3,-0.1)$) 	{{\color{red} $-2$}};
			\node at ($(g3.west) + (0.93,0.86)$) 	{{\color{red} $1$}};
			\node at ($(g4.east) + (0.3,0)$) 	{$\alpha_{4}$};
			\node at ($(g4.west) + (-0.35,-0.1)$) 	{{\color{red} $-2$}};
			\node at ($(g4.west) + (-0.79,0.86)$) 	{{\color{red} $1$}};		
			\node at ($(g5.east) + (0.3,0)$) 	{$\alpha_{5}$};
			\node at ($(g5.north) + (0,0.3)$) 	{{\color{red} $-2$}};
			\node at ($(g5.west) + (-0.95,-0.01)$) 	{{\color{red} $1$}};			
			\end{tikzpicture}
\caption{Dynkin diagram}
 \label{fig:Dynkin}
\end{center}
 {\footnotesize Numbers beside edges denote the intersection number $\langle \alpha_i, \check{\alpha}_j \rangle = \langle \alpha_j, \check{\alpha}_i \rangle$, while
numbers beside vertices denote the self-intersection number
$\langle \alpha_i, \check{\alpha}_i \rangle$.}
\end{figure}

Then, $w_{\alpha_i}$, $i=1, 2, 3, 4, 5$, acts on the N\'eron-Severi bi-lattice as
\begin{align*}
w_{\alpha_i}(D)= D-2 \frac{\langle D, \check{\alpha}_i   \rangle}{\langle \alpha_i, \check{\alpha}_i   \rangle} \alpha_i,\quad 
w_{\alpha_i}(d)=d-2\frac{\langle \alpha_i, d   \rangle}{\langle \alpha_i, \check{\alpha}_i   \rangle} \check{\alpha}_i
\end{align*}
for $D\in H^2({\cX_{\bf a}},\Z)$ and $d\in H_2({\cX_{\bf a}},\Z)$.

Thus, $w_{\alpha_i}$, $i=1, 2, 3, 4, 5$,  coincide with $w_{\kappa_0}, w_{\kappa_1}, w_{\kappa_\infty}, w_{\theta_1} ,w_{\theta_2}$ respectively, while $\sigma_j$, $j=1, 2, 3, 4$, act on the roots as transposition $(1,2)$, $(2,3)$, $(3,4)$, $(4,5)$ respectively.

\begin{rem}
The coroots are taken orthogonal to the vertical leaves, while there is not known automatic way to determine the roots, which are determined in a manner consistent with the action of the B\"acklund transformations on the bi-lattice.
\end{rem}

According to this remark, there seems to be no reasonable way to determine $\alpha_0$. Moreover, $w_{\alpha_0}$ was not a pseudo-isomorphism with 10 blow-ups.
In fact, if we apply the above formula to $w_{\alpha_0}$ with $D=\cH_q$, we obtain
$$\cH_q\mapsto  \cH_q+\frac{2}{5}(\cH_q+2\cH_r-\cE_1-\cE_3-\cE_5),$$ 
which is not in $H^2(\cX_\alpha,\Z)$, and $w_{\alpha_0}$ can not be realized as a birational map. 

Surprisingly, however, the following theorem holds and gives the reason why we take $\alpha_0$ as above.  

\begin{thm}
Kac's translation  (\S6.5 in \cite{Kac1990}):
$$T_{\alpha_i}(D) =D+\langle D,\check{\delta}\rangle\alpha_i +
\langle D,\check{\delta} - \check{\alpha}_i \rangle \delta $$
with\footnote{The author does not know geometric interpretation for those null roots.}
\begin{align}
\delta &= 2\alpha_0 +\sum_{i=1}^5 \alpha_i=3\cH_q+2\cH_r -\cE_{1,2,3,4,5,6,7,8,9,10}\\
\check{\delta} &= 2\check{\alpha}_0 +\sum_{i=1}^5 \check{\alpha}_i=
 2h_q+2h_r -e_{1,2,3,4,5,6,7,8,9,10}
\end{align}
can be realized as a pseudo-isomorphism for $i=1, 2, 3, 4, 5$. Moreover, $T_{\alpha_0}^2$ can be realized as a pseudo-isomorphism.
\end{thm}

Especially, $T_{\alpha_1}$ acts on $(\alpha_0,\alpha_1, \dots,\alpha_5)$ as
\begin{align*}
T_{\alpha_1}: & (\alpha_0,\alpha_1, \dots,\alpha_5) \mapsto
(\alpha_0,\alpha_1, \dots,\alpha_5) + \delta (-1,2,0,0,0,0)
\end{align*}
and 
 $T_{\alpha_0}^2$ acts on $(\alpha_0,\alpha_1, \dots,\alpha_5)$ as
\begin{align*}
T_{\alpha_0}^2:& (\alpha_0,\alpha_1, \dots,\alpha_5) \mapsto
(\alpha_0,\alpha_1, \dots,\alpha_5) + \delta (5,-2,-2,-2,-2,-2).
\end{align*}

\begin{proof}
Using Table~\ref{tab:pic}, it tuns out that the action of 
$$T_1 =\left(w_{\theta_2} w_{\theta_1} w_{\kappa_2} w_{\kappa_0}  w_{\alpha_0} \right)^2 $$
on the space with 21 blow-ups is trivial on the sub-lattice expanded by $\cE_{11}, \dots,  \cE_{21}$, and hence $T_1$ is a pseudo-isomorphism on the space with the first 10 blow-ups. Moreover, the action of $T_1$ on $\{\cX_{\alpha}\}$ coincides with $T_{\alpha_1}$.
 In other words, $T_{\alpha_1}$ is realized as $T_1$ as a birational map.

Obviously, we have
\begin{align*}
&T_{\alpha_1}=T_1,\quad T_{\alpha_2}=\sigma_1 T_1 \sigma_1,\quad  T_{\alpha_3}=\sigma_2\sigma_1 T_1 \sigma_1\sigma_2,\\ 
&T_{\alpha_4}=\sigma_3\sigma_2\sigma_1 T_1 \sigma_1 \sigma_2 \sigma_3,\quad T_{\alpha_5}=\sigma_4\sigma_3\sigma_2\sigma_1 T_1 \sigma_1\sigma_2\sigma_3\sigma_4,
\end{align*}
and using these translations,  we can realize $T_{\alpha_0}^2$ as
 $$
 T_{\alpha_0}^2= T_{-\alpha_1}T_{-\alpha_2}T_{-\alpha_3}T_{-\alpha_4}T_{-\alpha_5},
 $$ 
 where $T_{-\alpha_i}=T_{\alpha_i}^{-1}$ for $i=1,2,3,4,5$.
\end{proof}

\subsection*{Acknowledgments} About 20 years ago, the author had a study session with Hidetaka Sakai and Teruhisa Tsuda on the initial value space of the Garnier system based on Kimura's paper. Although not directly related to the present results, the author is grateful for the opportunity to learn Kimura's results and to have had a fruitful discussion.
The author also thanks the reviewer for careful reading of the manuscript. 
This work was supported by the Research Institute for Mathematical Sciences, an International Joint Usage/Research Center located in Kyoto University.

\end{document}